\documentclass[12pt]{article}
\usepackage{amsmath}
\usepackage{amssymb}

\newcommand{\fhat}{\hat{f}}
\newcommand{\fhatn}{\hat{f}(n)}
\newcommand{\ftilde}{{\tilde f}}
\newcommand{\Ftilde}{{\tilde F}}

\newcommand{\alext}{{\cal A}_c(\T)}
\newcommand{\balext}{{\cal B}_c(\T)}
\newcommand{\bvt}{{\cal BV}(\T)}

\newcommand{\DR}{{\cal D}(\R)}
\newcommand{\DPR}{{\cal D}'(\R)}
\newcommand{\bv}{{\cal BV}}
\newcommand{\intinf}{\int^\infty_{-\infty}}
\newcommand{\intpi}{\int_{-\pi}^{\pi}}

\newcommand{\N}{{\mathbb N}}
\newcommand{\R}{{\mathbb R}}
\newcommand{\C}{{\mathbb C}}
\newcommand{\Z}{{\mathbb Z}}
\newcommand{\T}{{\mathbb T}}
\newcommand{\fn}{\!:\!}

\newcommand{\qed}{\mbox{$\quad\blacksquare$}}
\newtheorem{theorem}{Theorem}
\newtheorem{lemma}[theorem]{Lemma}
\newtheorem{prop}[theorem]{Proposition}
\newtheorem{corollary}[theorem]{Corollary}

\newtheorem{defn}[theorem]{Definition}
\newtheorem{example}[theorem]{Example}

\begin{document}
\hspace{-2cm}
\raisebox{12ex}[1ex]{\fbox{{\footnotesize
Preprint
May 27, 2011.\quad
To appear in {\it Journal of Fourier Analysis and Applications}.
}}}

\begin{center}
{\large\bf Fourier series with the continuous primitive integral}
\vskip.25in
Erik Talvila\footnote{Supported by the
Natural Sciences and Engineering Research Council of Canada.
}\\ [2mm]
{\footnotesize
Department of Mathematics and Statistics \\
University of the Fraser Valley\\
Abbotsford, BC Canada V2S 7M8\\
Erik.Talvila@ufv.ca}
\end{center}

{\footnotesize
\noindent
{\bf Abstract.} 
Fourier series are considered on the one-dimensional torus 
for the space of periodic distributions
that are the distributional derivative of a continuous function.
This space of distributions is denoted $\alext$ and is a Banach
space under the Alexiewicz norm, 
$\|f\|_\T
=\sup_{|I|\leq 2\pi}|\int_I f|$, the supremum being taken
over intervals of length not exceeding $2\pi$.
It contains the periodic functions
integrable in the sense of Lebesgue and Henstock--Kurzweil.
Many of the properties of $L^1$ Fourier series continue to hold for
this larger  space, with the $L^1$ norm replaced by the Alexiewicz norm.
The Riemann--Lebesgue lemma takes the
form $\fhat(n)=o(n)$ as $|n|\to\infty$. The
convolution
is defined for $f\in\alext$ and $g$ a periodic function of bounded
variation.  The convolution
commutes with translations and
is commutative and associative.  There is the estimate
$\|f\ast g\|_\infty\leq \|f\|_\T \|g\|_\bv$. 
For $g\in L^1(\T)$,
$\|f\ast g\|_\T\leq \|f\|_\T \|g\|_1$.
As well, 
$\widehat{f\ast g}(n)=\fhatn \hat{g}(n)$. 
There are versions of the
Salem--Zygmund--Rudin--Cohen factorization theorem, Fej\'er's lemma and
the Parseval equality.
The trigonometric polynomials
are dense in $\alext$.
The convolution of $f$ with 
a sequence of summability kernels converges to $f$ in the Alexiewicz norm.
Let $D_n$ be the Dirichlet kernel and let
$f\in L^1(\T)$.  Then $\|D_n\ast f-f\|_\T\to 0$ as $n\to\infty$.
Fourier
coefficients of functions of bounded variation are characterized.  An appendix contains a type of
Fubini theorem.
\\
{\bf 2000 subject classification:} 26A39, 42A16, 46F10
}

\section{Introduction and notation}\label{introduction}
In this paper we consider Fourier series on the one dimensional
torus.  Progress in Fourier analysis has gone hand in hand with progress
in theories of integration.  This is perhaps best exemplified by the work
of Riemann and Lebesgue using the integrals named after them.
We describe
below the {\it continuous primitive integral}.  This is an integral
that includes the Lebesgue, Henstock--Kurzweil and wide Denjoy integrals.
It has a simple definition in terms of distributions.  The space of
distributions integrable in this sense is a Banach space under the
Alexiewicz norm.  Many properties of Fourier series that hold for $L^1$
functions continue to hold in this larger space with the $L^1$ norm
replaced by the Alexiewicz norm.

We use the following notation for distributions.  
The space of {\it test functions} is
$\DR=C^\infty_c(\R)=\{\phi\fn\R\to\R\mid \phi\in C^\infty(\R) \text{ and }
{\rm supp}(\phi) \text{ is compact}\}$.  The {\it support} of 
function $\phi$ is the
closure of the set on which $\phi$ does not vanish and is denoted
${\rm supp}(\phi)$.  Under usual pointwise operations $\DR$ is 
a linear space over field $\R$.  In $\DR$ we have
a notion of convergence.  If $\{\phi_n\}\subset\DR$ then $\phi_n\to 0$ as
$n\to\infty$
if there is a compact set $K\subset\R$ such that for each $n$,
${\rm supp}(\phi_n)\subset K$, and for each $m\geq 0$ we have
$\phi_n^{(m)}\to 0$ uniformly on $K$ as $n\to\infty$.  The {\it distributions}
are denoted $\DPR$ and
are the continuous linear functionals on $\DR$.  For $T\in\DPR$ and
$\phi\in\DR$ we write $\langle T,\phi\rangle\in\R$.  For $\phi,\psi\in\DR$
and $a,b\in\R$ we have $\langle T,a\phi+b\psi\rangle=a\langle T,\phi\rangle
+ b\langle T,\psi\rangle$.  And, if $\phi_n\to 0$ in $\DR$ then
$\langle T, \phi_n\rangle\to 0$ in $\R$.  Linear operations are defined
in $\DPR$ by $\langle aS+bT,\phi\rangle=a\langle S,\phi\rangle
+ b\langle T,\phi\rangle$ for $S,T\in\DPR$; $a,b\in\R$ and $\phi\in\DR$.
If $f\in L^1_{\it loc}$
then $\langle T_f, \phi\rangle = \intinf f(x)\phi(x)\,dx$ defines a distribution
$T_f\in\DPR$.  The integral exists as a Lebesgue integral.
All distributions have derivatives of
all orders that are themselves distributions.
For $T\in\DPR$ and
$\phi\in\DR$ the distributional derivative of $T$ is $T'$
where $\langle T',\phi\rangle = - \langle T,\phi'\rangle$.
If $p\fn\R\to\R$ is a function that is differentiable in
the pointwise sense at $x\in\R$ then we write its derivative
as $p'(x)$.
If $p$ is
a $C^\infty$ bijection such that
$p'(x)\not=0$ for any $x\in\R$  then the composition with
distribution $T$ is defined by
$\langle T\circ p,\phi\rangle =
\langle T,\frac{\phi\circ p^{-1}}{p'\circ p^{-1}}\rangle$ for all
$\phi\in\DR$. 
Translations are a special case.
For $x\in\R$ define the {\it translation} $\tau_x$ on distribution
$T\in\DPR$ by $\langle \tau_xT, \phi\rangle = \langle T,\tau_{-x}\phi\rangle$
for test function $\phi\in\DR$ where $\tau_x\phi(y)=\phi(y-x)$.
A distribution $T\in\DPR$ is {\it periodic} if 
$\langle \tau_pT,\phi\rangle=\langle T,\phi\rangle$ for
some $p>0$ and all $\phi\in\DR$.  The least such positive $p$ is
the {\it period}.  In this paper, periodic will
always mean periodic with period $2\pi$.
Periodic distributions are defined in an alternative manner
in \cite{katznelson} and \cite{zemanian}.
All of the results on distributions we use can be
found in these works and \cite{friedlander}.

We define the torus $\T=\{z\in\C\mid |z|=1\}=
\{e^{i\theta}\mid\theta\in\R\}$.  The real interval $[-\pi,\pi)$ will
be used as a model for $\T$. 

The continuous primitive integral was discussed on the real line in
\cite{talviladenjoy}.  As the name suggests, this integral is characterized
by having a primitive that is a continuous function;  the integrable
distributions are those that are the distributional derivative
of a continuous function.
Take the space of primitives as
$\balext=\{F\fn\R\to\R|F\in C^0(\R), F(-\pi)=0,
F(x)=F(y)+nF(\pi)
\text{ if } y\in[-\pi,\pi), x=y+2n\pi \text{ for } n\in\Z\}$. 
Note that $F\in\balext$ is periodic on $\R$ if and only if $F(\pi)=0$.
If $x\in\R$ and $n\in\Z$ then $F(x+2n\pi)=F(x)+nF(\pi)$ and
$F(x)=(x-x\,{\rm mod}\ 2\pi)F(\pi)/(2\pi) +F(x\,{\rm mod}\ 2\pi)$.
It is easy to see that $\balext$ is a Banach space under the
uniform norm $\|F\|_{\T,\infty}=
\sup_{|\alpha-\beta|\leq 2\pi}|F(\alpha)-F(\beta)|$.
The integrable distributions on the torus are then given by
$\alext=\{f\in\DPR\mid f=F' \text{ for some } F\in\balext\}$.  
For $a,b\in\R$ the integral
of $f\in\alext$ is $\int_a^bf=F(b)-F(a)$ where $F\in\balext$ and
$F'=f$.  Note that for all $a,b\in\R$ and all
$m,n\in\Z$ we have $\int_{a+2m\pi}^{b+2n\pi} f=\int_a^bf + (n-m)\int_{-\pi}^{\pi}f$.  If $f$ is complex-valued, the real and imaginary parts
are integrated separately.  The distributional
differential equation $T'=0$ has only constant solutions and
we have made our primitives in $\balext$ vanish at $-\pi$ so the primitive
of a distribution in $\alext$ is unique.  

If $f\fn\R\to\R$ is
a periodic function that is locally integrable in the
Lebesgue, Henstock--Kurzweil or wide Denjoy sense then $T_f\in\alext$.
Thus, if $f(t) = t^{-2}\cos(t^{-2})$ for $t\in(0,\pi)$ and $f(t)=0$
for $t\in[-\pi,0]$ with $f$ extended periodically, then $T_f\in\alext$
but $f$ is not Lebesgue integrable.  In this case, $f$ has an
improper Riemann integral and $0$ is the only point of nonabsolute
summability.  There are examples of functions integrable in the
Henstock--Kurzweil sense but not in the Lebesgue sense for which
the set of points of nonabsolute
summability has positive measure.  See \cite{jeffery}.  
Such functions correspond to
distributions 
integrable in the  continuous primitive sense.
We will usually drop the distinction between $f$ and $T_f$.
As well, if $F\in\balext$ is a function of Weierstrass type that is
continuous but has a pointwise derivative nowhere then the distributional
derivative of $F$ exists and $F'\in\alext$.  If $F$ is a continuous
singular function, so that $F'(x)=0$ a.e., then $F'\in\alext$ and
the continuous primitive integral is $\int_a^bF'=
F(b)-F(a)$.  In this case, $F'\in L^1(\T)$ but the Lebesgue
integral gives $\int_a^bF'(x)\,dx=0$.

If $f\in\alext$
and $F\in\balext$ is its primitive then the action of
$f$ on test function $\phi\in\DR$ is given by $\langle f,\phi\rangle
=\langle F',\phi\rangle=-\langle F,\phi'\rangle=-\intinf F(x)\phi'(x)\,dx$.
This last integral exists as a Riemann integral.
And,  for $f\in\alext$ with primitive $F\in\balext$,
\begin{align*}
&\langle \tau_{2\pi}f,\phi\rangle=\langle f,\tau_{-2\pi}\phi\rangle
=-\intinf F(x)\phi'(x+2\pi)\,dx\\
&=-\intinf F(x-2\pi)\phi'(x)\,dx
=-\intinf F(x)\phi'(x)\,dx+F(\pi)\intinf\phi'(x)\,dx\\
&=-\langle F,\phi'\rangle=\langle f,\phi\rangle,
\end{align*}
so $f$ is periodic.  If $F\in C^0(\T)$ then $F'\in\alext$.
Note that distributions in $\alext$ are
tempered and of order one.  See \cite{friedlander} for the
definitions.

Distributions in $\alext$ can be composed with continuous
functions and this leads to a very powerful change of variables
formula.  See \cite[Theorem~11]{talviladenjoy}.

The {\it Alexiewicz norm} of $f\in\alext$ is $\|f\|_\T
=\sup_{|I|\leq 2\pi}|\int_I f|$, the supremum being taken
over intervals of length not exceeding $2\pi$.  We have
$\|f\|_\T=\|F\|_{\T,\infty}
=\max_{|\beta-\alpha|\leq 2\pi}|F(\beta)-F(\alpha)|$
where $F\in\balext$ is the primitive of $f$.  The integral provides
a linear isometry and isomorphism between $\alext$ and $\balext$.  Define
$\Phi\fn\alext\to\balext$ by $\Phi[f](x)=\int_{-\pi}^xf$.  Then
$\Phi$ is a linear bijection and $\|f\|_\T=\|\Phi[f]\|_{\T,\infty}$.
Hence, $\alext$
is a Banach space.  The spaces of periodic Lebesgue, Henstock--Kurzweil
and
wide Denjoy integrable functions are not complete under the
Alexiewicz norm.  The space $\alext$ furnishes their completion.
An equivalent norm is $\|f\|'_\T=\sup_{-\pi\leq x\leq \pi}|\int_{-\pi}^xf|$.

The multipliers and dual space of $\alext$ are given by the functions
of bounded variation.  If $g\fn\R\to\R$ is periodic then its variation
over $\T$ is given by $Vg=\sup\sum|g(s_i)-g(t_i)|$ where the supremum
is taken over all disjoint intervals $\{(s_i,t_i)\}\subset(-\pi,\pi)$.
We write $\bvt$ for the periodic functions with finite variation.
This is a Banach space under the norm $\|g\|_\bv=\|g\|_\infty+Vg$.
If $g$ is complex-valued with real and imaginary
parts $g_r$ and $g_i$, then $Vg=\sqrt{(Vg_r)^2+(Vg_i)^2}$.
If $f\in\alext$ with primitive $F\in\balext$ and $g\in\bvt$ then the integral of $fg$ is defined
using a Riemann--Stieltjes integral, 
\begin{equation}
H(x)=\int_{-\pi}^x fg=F(x)g(x)-\int_{-\pi}^xF(t)\,dg(t),
\quad x\in[-\pi,\pi).\label{parts}
\end{equation}  Extension of $H$ outside this interval using
$H(x)=(x-x\,{\rm mod}\ 2\pi)H(\pi)/(2\pi) +H(x\,{\rm mod}\ 2\pi)$
yields an element of $\balext$ whose derivative is then interpreted
as $fg\in\alext$.  Note that $\bvt\subset L^1(\T)\subset\alext$.

Growth estimates and other basic properties of Fourier coefficients are
proved in
Theorem~\ref{coefficients}. Let $\fhat(n)=\intpi f(t) e^{-int}\,dt$ 
denote the Fourier coefficient
of $f\in\alext$.  The Riemann--Lebesgue lemma takes the
form $\fhat(n)=o(n)$ as $|n|\to\infty$.  In Theorem~\ref{bvconv}, the
convolution $f\ast g(x)=\intpi f(x-t)g(t)\,dt$
is defined for $f\in\alext$ and $g\in\bvt$.  The convolution is
then continuous,
commutes with translations and
is commutative and associative.  There is the estimate
$\|f\ast g\|_\infty\leq \|f\|_\T \|g\|_\bv$.  As well, 
$\widehat{f\ast g}(n)=\fhatn \hat{g}(n)$.  The integral
$\intpi f(x-t)g(t)\,dt$ need not exist for $f\in\alext$ and $g\in L^1(\T)$.
But using the density of $L^1(\T)$ in $\alext$ and  the density
of $\bvt$ in $L^1(\T)$ we can define the convolution for $f\in\alext$
and $g\in L^1(\T)$ as the limit of a sequence of convolutions $f_k\ast g$
for $f_k\in L^1(\T)$ or as the limit of $f\ast g_k$ for
$g_k\in \bvt$ (Theorem~\ref{L1conv}).  The usual properties of convolution
continue to hold.  Now we have $\|f\ast g\|_\T\leq \|f\|_\T \|g\|_1$.
Theorem~\ref{factorization} gives a version of the
Salem--Zygmund--Rudin--Cohen factorization theorem, $\alext=L^1(\T)\ast \alext$.
Using the  Fej\'er kernel it is shown that the trigonometric polynomials
are dense in $\alext$.  There is the
uniqueness result that if $\fhat={\hat g}$
then $f=g$ as distributions in $\alext$.
As well,
the convolution of $f$ with 
a sequence of Fej\'er kernels converges to $f$ in the Alexiewicz norm
(Theorem~\ref{uniqueness}).
Let $D_n$ be the Dirichlet kernel and let
$f\in L^1(\T)$.  Then $\|D_n\ast f-f\|_\T\to 0$ as $n\to\infty$ (Theorem~
\ref{dirichlet}).
Example~\ref{Dnexample} shows there is
$f\in\alext$ such that $\|D_n\ast f-f\|_\T\not\to 0$.
Proposition~\ref{fejerlemma} is a version of Fej\'er's lemma
and Theorem~\ref{parseval} is a type of Parseval equality.
Theorem~\ref{bvcharacterization} gives a characterization of Fourier
coefficients of functions in $\bvt$.  An appendix contains a type of
Fubini theorem.

We will use the following
version of the H\"older inequality from
the Appendix of
\cite{talvilafourier}.
\begin{prop}[H\"older inequality]\label{holder}
Let $f\in\alext$. If $g\in\bvt$ then
$\left|\int_{-\pi}^\pi fg\right|  \leq |\int_{-\pi}^\pi f|\inf|g|+\|f\|_\T Vg
\leq  
\|f\|_\T\|g\|_\bv$.
\end{prop}

Distributions in $\alext$ are continuous in the Alexiewicz norm.  This
means that if $f\in\alext$ then $\|f-\tau_sf\|_\T\to 0$ as $s\to 0$.
See \cite[Theorem~28]{talviladenjoy} for a proof.

A function on  the real line is called {\it regulated} if it has a 
left limit and a right limit at each point. 
The {\it regulated primitive integral} integrates those distributions
that are the distributional derivative of a regulated function.
Analogous to $\alext$, the
space of integrable distributions is a Banach space.  This space
includes $\alext$ and also all signed Radon measures.  A theory
of Fourier series can be obtained as in the present paper. The
chief difference between these two integrals is in
the integration by parts formula and in the fact that we no longer
have continuity in the Alexiewicz norm.  See \cite{talvilarpi}.

If $u$ is a periodic distribution then it has a Fourier series
given by
$u(x)=[1/(2\pi)]\sum_{n\in\Z}\hat{u}_ne^{inx}$.  The convergence is in
the distributional sense, 
i.e., weak convergence.
There is also a converse, any trigonometric series with coefficients
of polynomial growth is the Fourier series of a distribution. 
See \cite[Theorem~8.5.2, Theorem~8.5.3]{friedlander}.  Other methods
of defining Fourier series of distributions are given in 
\cite{edwardsII} and \cite{zemanian}.  Since $\alext$ is a subspace
of distributions, all of the results in these works continue to hold.
However, $\alext$ is also a Banach space.  We will see below that
Fourier series of distributions in $\alext$ behave more like those for
$L^1$ functions than for general distributions.  

In this paper we develop basic properties of Fourier series and
convolutions ab initio from the definition
of the integral.  As the functions of bounded variation are pointwise
multipliers for distributions in $\alext$ it follows that $\alext$ is
a Banach $\bvt$-module over the pointwise algebra of $\bvt$.
And, as is shown in Section~\ref{sectionconvolution} below, distributions in
$\alext$ can be convolved with functions in $L^1(\T)$ such that $\alext$ is
a Banach $L^1(\T)$-module over the convolution algebra of
$L^1(\T)$.  Although we have employed a concrete approach here, such a 
two-module property may allow the abstract methods
developed in \cite{braunfeichtinger} to be used to deduce some of the
theorems below.  An anonymous referee suggested Corollary~\ref{bvseries}
might be proved this way.

\section{Fourier coefficients}
Let $e_n(t)=e^{int}$.
If $f\in\alext$ then the Fourier coefficients of $f$ are
$\fhatn=\langle f,e_{-n}\rangle=\intpi f e_{-n}=\intpi f(t)e^{-int}\,dt$, where $n\in\Z$.  
Since the functions
$e_n$ and $1/e_{n}$ are in $\bvt$ for each $n\in\Z$, the
Fourier coefficients exist on $\Z$ as continuous primitive integrals if and
only if $f\in\alext$.  Let $F(x)=\int_{-\pi}^x f$
be the primitive of $f$.  Integrating by parts as in  \eqref{parts}
gives 
\begin{equation}
\fhatn=(-1)^nF(\pi)+in\intpi F(t)e^{-int}\,dt.\label{fnparts}
\end{equation}
This last integral is
the Riemann integral of a continuous function.  Formula \eqref{fnparts}
can be used as an
alternative definition of $\fhatn$.  Note also that
$\fhatn=\int_{\alpha}^{\alpha+2\pi} f(t)e^{-int}\,dt$ for each $\alpha\in\R$.  
The following properties
of the Fourier coefficients follow easily from the linearity of
the integral and from \eqref{fnparts}.  The complex conjugate is
denoted $\overline{x+iy}=x-iy$ for $x,y\in\R$.  We will  take  $f$ to
be real-valued but only trivial changes are required for complex-valued
distributions.
\begin{theorem}\label{coefficients}
Let $f, g\in\alext$.  Then (a) $\widehat{f+g}(n)=\fhatn+\hat{g}(n)$;
(b) if $\alpha\in\C$ then $\widehat{(\alpha f)}(n)=\alpha\fhatn$;
(c) $\overline{\fhat(n)} =\fhat(-n)$;
(d) if  $s\in\R$ then $\widehat{\tau_sf}(n)=\fhatn e^{-ins}$;
(e) $|\fhatn|\leq |F(\pi)|+ |n|\intpi|F|$ where $F(x)=\int_{-\pi}^xf$;
(f) for $n\not= 0$, $|\fhatn|\leq 4\sqrt{2}\,|n|\|f\|_\T$;
(g) $\fhatn=o(n)$ as $|n|\to\infty$ and this estimate is sharp;
(h) for $n\not=0$ we have $|\fhatn|\leq 2\sqrt{2}\,|n|\|f-\tau_{\pi/n}f\|_\T$;
(i) if $F\in C^0(\T)$ then $\widehat{F'}(n)=in\hat{F}(n)$;
(j) if $F\in C^{k-1}(\T)$ for some $k\in\N$ then for $n\not=0$ and
each $0\leq \ell\leq k$,
$\hat{F}(n)=(in)^{-\ell}\widehat{F^{(\ell)}}(n)$ and 
$$|\hat{F}(n)|\leq 4\sqrt{2}
\min_{0\leq \ell\leq k}\frac{\|F^{(\ell)}\|_\T}{|n|^{\ell-1}}.
$$
As $|n|\to\infty$, $\hat{F}(n)=o(n^{1-k})$.
\end{theorem}
Part (g) is a version of the Riemann--Lebesgue lemma for the continuous
primitive integral.
When $f\in L^1(\T)$ then $\fhatn=o(1)$ as $|n|\to\infty$.  This estimate
is sharp in the sense that if $\psi\fn\N\to(0,\infty)$ and
$\psi(n)=o(1)$ as $n\to\infty$ then there is a function $f\in L^1(\T)$
such that $\fhatn\not=o(\psi(n))$ as $|n|\to\infty$.  Estimates similar
to those in (j) appear in \cite[I~4.4]{katznelson} for $F^{(k-1)}$
absolutely continuous.

\bigskip
\noindent
{\bf Proof:}
To prove (f), apply the H\"older
inequality (Proposition~\ref{holder}) as follows. Notice that
the minimum of $|\sin(nt)|$ and $|\cos(nt)|$ for $|t|\leq \pi$
are both zero.  Then
$$|\fhatn|  \leq  \|f\|_\T\sqrt{
\left[\intpi|n\sin(nt)|\,dt\right]^2 +
\left[\intpi|n\cos(nt)|\,dt\right]^2}
  =  4\sqrt{2}\,|n|\|f\|_\T.
$$
Part (g) 
follows upon integrating by parts and using the $L^1$ form of the
Riemann--Lebesgue
lemma on the integral  $\intpi F(t)e^{-int}\,dt$.  The estimate
was proved sharp by Titchmarsh \cite{titchmarsh}.
To prove (h), use a linear change of variables to write
$\fhatn=(1/2)\intpi[f(t)-f(t-\pi/n)]e^{-int}\,dt$. 
Then proceed as in (f) to get $|\fhatn|\leq 2\sqrt{2}\,|n|\|f-\tau_{\pi/n}f\|_\T$.
Since $f$ is continuous in the Alexiewicz norm this also gives the
little oh estimate in (g).  Part (i) follows from integrating by parts
and then (j) is obtained using (i) with the estimates in (f) and (g). \qed

The quantity $\omega(f,\delta)=\sup_{|t|<\delta}\|f-\tau_tf\|_\T$
is known as the {\it modulus of continuity} of $f$ in the Alexiewicz
norm.  Part (h) gives $|\fhatn|\leq 2\sqrt{2}\,|n|\omega(f,\pi/|n|)$.

Katznelson \cite{katznelson} gives other estimates for $\fhatn$
under such assumptions as $f$ is of bounded variation, absolutely
continuous, Lipschitz continuous or  in $L^p(\T)$. 

The next theorem shows that when we have a sequence converging in
the Alexiewicz norm, the Fourier coefficients also converge.
\begin{theorem}\label{coefficient convergence}
For $j\in\N$, 
let $f, f_j\in\alext$ such that $\|f_j-f\|_\T\to 0$ as $j\to\infty$.
Then for each $n\in\Z$  we have 
$\hat{f_j}(n)\to \fhatn$ as $j\to\infty$.
The convergence need not be uniform in $n\in\Z$.
\end{theorem}
\bigskip
\noindent
{\bf Proof:}
If $n=0$ then $|\hat{f_j}(0) -\fhat(0)|=\left|\intpi\left[f_j(t)-f(t)\right]
\,dt\right|\leq \|f_j-f\|_\T$.  If $n\not= 0$ then
\begin{eqnarray}
|\hat{f_j}(n) -\fhat(n)| & = & \left|\intpi [f_j(t)-f(t)]e^{-int}\,dt\right|
\notag\\
& \leq & 4\sqrt{2}\, |n|\|f_j-f\|_\T\to 0 
\text{ as } j\to\infty.\label{normestholder}
\end{eqnarray}
Theorem~\ref{coefficients}(f) is used in \eqref{normestholder}.

To show the convergence need not be uniform, let $f_j(t)=e^{ijt}$ and
$f=0$.  Then $\|f_j\|_\T=\sup_{|\alpha-\beta|\leq 2\pi}
|\int_{\alpha}^\beta e^{ijt}\,dt|$.
We have
$$
\left|\int_{\alpha}^\beta e^{ijt}\,dt\right|  =  \frac{1}{j}\left|
e^{ij(\beta-\alpha)}-1\right|
 \leq \frac{2}{j}.
$$
Equality is realized when $\beta=\alpha+\pi/j$.  Hence, $\|f_j\|_\T=2/j\to 0$
as $j\to\infty$.  But, $\hat{f_n}(n)=\intpi dt=2\pi\not\to 0$.\qed

This is different from the case $f, f_j\in L^1(\T)$.  There, if $\{f_j\}$
converges to $f$ in the $L^1$ norm then $\hat{f_j}(n)$ converges to $\fhatn$ 
uniformly in $n$ as $j\to\infty$.  See \cite[I~Corollary~1.5]{katznelson}.

\section{Convolution}\label{sectionconvolution}
The convolution is one of the most important operations in analysis,
with applications to differential equations, integral equations
and approximation of functions.  For $f\in\alext$ and $g\in\bvt$
the convolution is
$\intpi(f\circ r_x)g$ where $r_x(t)=x-t$.  We write this as
$f\ast g(x)=\intpi f(x-t)g(t)\,dt$.  This integral
exists for all such $f$ and $g$.  The convolution inherits smoothness
properties from $f$ and $g$.  We will also use a limiting process
to define the convolution for $g\in L^1(\T)$.  This then makes
$\alext$ into an $L^1(\T)$-module over the $L^1(\T)$ convolution algebra.
See \cite[32.14]{hewittrossII} for the definition.
The convolution was
considered for the continuous primitive integral on the real line
in \cite{talvilaconv}.   Many of the results of that paper are easily
adapted to the setting of $\T$, especially differentiation and integration
theorems which we do not reproduce here.

When $f$ and $g$ are in $L^1(\T)$ the convolution $f\ast g$ is
commutative and associative.  The estimate $\|f\ast g\|_1\leq
\|f\|_1\|g\|_1$ shows the convolution is a bounded linear operator
$\ast\fn L^1(\T)\times L^1(\T)\to L^1(\T)$ and $L^1(\T)$ is
a Banach algebra under convolution.  See \cite{katznelson} for
details.

Since $\bvt$ is the dual of $\alext$, many of the usual properties
of convolutions hold when it is defined on $\alext\times\bvt$.

\begin{theorem}\label{bvconv}
Let $f\in\alext$ and let $g\in\bvt$.  Then (a) $f\ast g\in C^0(\T)$
(b)
$f\ast g=g\ast f$
(c)
$\|f\ast g\|_\infty\leq \|f\|_\T \|g\|_\bv$
(d)
for $y\in\R$ we have $\tau_y(f\ast g)=(\tau_yf)\ast g
=f\ast(\tau_yg)$.
(e) If $h\in L^1(\T)$ then $f\ast (g\ast h)=(f\ast g)\ast h\in C^0(\T)$.
(f) For each $f\in\alext$, define $\Phi_f\fn\bvt\to C^0(\T)$ 
by $\Phi_f[g]=f\ast g$.
Then $\Phi_f$ is a bounded linear operator and $\|\Phi_f\|\leq\|f\|_\T$.
For each $g\in\bvt$, define $\Psi_g\fn\alext\to C^0(\T)$ 
by $\Psi_g[f]=f\ast g$.
Then $\Psi_f$ is a bounded linear operator and $\|\Psi_g\|\leq\|g\|_\bv$.
(g) We have $\widehat{f\ast g}(n)=\fhatn \hat{g}(n)$ for all $n\in\Z$.
(h) $\|f\ast g\|_\T\leq\|f\|_\T\|g\|_1$.
(i) Let $f,g\in L^1(\T)$.  Then $\|f\ast g\|_\T\leq \|f\|_\T\|g\|_1\leq
\|f\|_1\|g\|_1$.
\end{theorem}
\bigskip
\noindent
{\bf Proof:}
The proofs of (a) through (e) are essentially the same as for 
\cite[Theorem~1]{talvilaconv}. Proposition~\ref{holder} is used
in the proof of (c).  Theorem~\ref{fubini} is used in
the proof of (e).  Part (f) follows from  part (c).
To prove (g), write $\widehat{f\ast g}(n)=\intpi\intpi
f(x-t)e^{-in(x-t)}g(t)e^{-int}\,dt\,dx$.  By Theorem~\ref{fubini}
we can interchange the orders of integration.  The result then
follows upon a change of variables.  To prove (h), use Theorem~\ref{fubini}
and a linear change of variables
to write
\begin{eqnarray*}
\int_\alpha^\beta f\ast g(x)\,dx & = &  
\int_\alpha^\beta\intpi f(x-t)g(t)\,dtdx
=  \intpi\int_\alpha^\beta f(x-t)g(t)\,dxdt\\
 & = & \intpi g(t)\int_{\alpha-t}^{\beta-t} f(x)\,dxdt.
\end{eqnarray*}
Then 
$$
\left|\int_\alpha^\beta f\ast g(x)\,dx\right|\leq
\sup_{u<v}\left|\int_u^v f\right|\|g\|_1
\leq \|f\|_\T\|g\|_1.
$$
The proof of (i) is similar but now the usual Fubini theorem is used.\qed

Using two equivalent norms, we can have equality in part (f).
Define $\|f\|'_\T=\sup_{-\pi\leq x\leq \pi}|\int_{-\pi}^xf|$ for
$f\in\alext$ and define $\|g\|'_\bv=|g(-\pi)|+0.5Vg$ for $g\in\bvt$.
These norms are equivalent to $\|\cdot\|_\T$ and $\|\cdot\|_\bv$, respectively.
Given $f\in\alext$ with $f\not=0$ there is $\alpha\in(-\pi,\pi]$ such
that $\|f\|'_\T=|\int_{-\pi}^\alpha f|$.  Define $g\in\bvt$ by
$g(t)=\chi_{(-\alpha,\pi)}(t)$ for $t\in[-\pi,\pi)$ and extend
periodically.  Then $\|g\|'_\bv=1$ and $|f\ast g(0)|=\|f\|'_\T$.
With these norms, $\|\Phi_f\|=\|f\|'_\T$.   However, we can have
$\|\Psi_g\|<\|g\|_\bv$.  Let $g(t)=(1/3)\chi_{\{0\}}(t)$ for
$t\in[-\pi,\pi)$ and extend
periodically.  Then $\|g\|_\bv=1$ but $f\ast g=0$ for each $f\in \alext$.
Hence, $\|\Psi_g\|=0$.  This problem goes away if we replace $\bvt$ with
functions of normalised bounded variation.  Fix $0\leq\lambda\leq 1$.
A function $g\in\bvt$ is of normalised bounded variation if
$g(x)=(1-\lambda)g(x-)+\lambda g(x+)$ for each $x\in[-\pi,\pi)$.  The
case $\lambda=0$ corresponds to left continuity and $\lambda=1$ 
corresponds to right continuity.  See \cite{talvilarpi} for
details.  If $g$ is then normalised for some $\lambda$ such that
$g(-\pi)=0$ then $\|g\|'_\bv=Vg$.  Similar to the method in the proof of 
Theorem~\ref{L1conv}(e) it can be seen that $\|\Psi_g\|=Vg$.

Linearity in each component, associativity (e) and inequality (c) 
show that $\alext$ 
is a $\bvt$-module.  Note that $\bvt$ is a Banach algebra under
pointwise operations.

\begin{example}
{\rm
Note that $f\ast g$ need not be of bounded variation and hence
need not be absolutely continuous.  For example, let $g$ be
the periodic extension of $\chi_{(0,\pi)}$.  Then
$f\ast g(x)=F(x)-F(x-\pi)$ where $F\in\balext$ is the primitive of $f$.
Since $F$ need not be of bounded variation, the same can be said for
$f\ast g$.  For instance, take $F(x)=x\sin(x^{-2})\chi_{[0,\pi)}(x)$ on $[-\pi,\pi)$
with $F(0)=0$.
}
\end{example}
Using the estimate $\|f \ast g\|_\T\leq \|f\|_\T \|g\|_1$ 
(Theorem~\ref{bvconv}(h)) and the
fact that $\bvt$ is dense in $L^1(\T)$ we can define
$f\ast g$ for $f\in\alext$ and $g\in L^1(\T)$.
Since $L^1(\T)$ is dense in $\alext$ we can also define the
convolution using a sequence in $L^1(\T)$ with the inequality
in Theorem~\ref{bvconv}(i).
\begin{defn}\label{abdefn}
Let $f\in\alext$ and $g\in L^1(\T)$. (a)  Let $\{g_k\}\subset\bvt$
such that $\|g_k-g\|_1\to 0$.  Then $f\ast g$ is the unique element of
$\alext$ such that $\|f\ast g_k-f\ast g\|_\T\to 0$.
(b) Let $\{f_k\}\subset L^1(\T)$
such that $\|f_k-f\|_\T\to 0$.  Then $f\ast g$ is the unique element of
$\alext$ such that $\|f_k\ast g-f\ast g\|_\T\to 0$.
\end{defn}
See \cite{talvilaconv} for a proof that (a) defines a unique element
of $\alext$.
The validity of (b) and the
equality of definitions (a) and (b) is proved
following  the proof of Theorem~\ref{L1conv}.

\begin{theorem}\label{L1conv}
Let $f\in\alext$ and $g\in L^1(\T)$.  Then (a) $f\ast g\in \alext$
(b) $\|f\ast g\|_\T\leq \|f\|_\T \|g\|_1$ 
(c) for $y\in\R$ we have $\tau_y(f\ast g)=(\tau_yf)\ast g
=f\ast(\tau_yg)$.
(d)  If $h\in L^1(\T)$ then $f\ast (g\ast h)=(f\ast g)\ast h\in \alext$.
(e) For each $f\in\alext$, define
$\Phi_f\fn L^1(\T)\to \alext$ by $\Phi_f[g]=f\ast g$.
Then $\Phi_f$ is a bounded linear operator and $\|\Phi_f\|=\|f\|_\T$.
For each $g\in L^1(\T)$, define
$\Psi_g\fn \alext\to \alext$ by $\Psi_g[f]=f\ast g$.
Then $\Psi_g$ is a bounded linear operator and $\|\Psi_g\|=\|g\|_1$.
(f) We have $\widehat{f\ast g}(n)=\fhatn \hat{g}(n)$ for all $n\in\Z$.
\end{theorem}
\bigskip
\noindent
{\bf Proof:}
Using Theorem~\ref{fubini}, 
the proofs of (a) through (d) are essentially the same as for 
\cite[Theorem~3]{talvilaconv}, taking Theorem~\ref{bvconv} into account.
To proof (e), let $f\in\alext$.  Then
$\|\Phi_f\|=\sup_{\|g\|_1=1}\|f\ast g\|_\T\leq\sup_{\|g\|_1=1}\|f\|_\T\|g\|_1
=\|f\|_\T$.  To show we have equality, let $k_n\in\bvt$ be a sequence of
non-negative summability
kernels as in Definition~\ref{defnsummability}.  
Using Theorem~\ref{normconv}, $\|\Phi_f\|\geq\|f\ast k_n\|_\T
\to\|f\|_\T$ as $n\to\infty$.
The inequality $\|\Psi_g\|\leq \|g\|_1$ follows
similarly.  To prove we have equality, note that $C^0(\T)$ is dense in
$L^1(\T)$ and such functions are uniformly continuous, so for each
$g\in L^1(\T)$ and each
$\epsilon>0$ there is a step function $\sigma(x)=\sum_{j=1}^n\sigma_j
\chi_{(a_{j-1},a_j)}(x)$ defined by a uniform partition 
$-\pi=a_0<a_1<\cdots<a_n=\pi$, $a_j=-\pi+2\pi j/n$, for which $\|g-\sigma\|_1
<\epsilon$.  It then suffices to find a sequence $\{f_m\}\subset\alext$
with $\|f_m\|_\T=1$ and $\|f_m\ast\sigma\|_\T\to\|\sigma\|_1$ as $m\to\infty$.
Define $f_m(x)=\sum_{i=1}^n\epsilon_i p_m(x+a_{i-1})$ where
$\epsilon_i={\rm sgn}(\sigma_i)$ and $p_m$ is a delta sequence.  This is a
sequence of continuous  functions
$p_m\geq 0$ such that there is a sequence of real numbers $\delta_m\downarrow 0$
with ${\rm supp}(p_m)\subset(-\delta_m,\delta_m)$ and 
$\int_{-\delta_m}^{\delta_m}p_m(x)\,dx=1$.  Take $m$ large enough so that
$\delta_m<\pi/n$.  With the Fubini theorem we have
\begin{equation}
\|f_m\ast\sigma\|_\T  \geq  \int_0^{2\pi/n}f_m\ast\sigma(x)\,dx
=\sum_{i,j=1}^n\epsilon_i\sigma_j\int_{a_{j-1}}^{a_j}\int_{a_{i-1}-t}^{a_{i}
-t}p_m(x)\,dx\,dt.\label{pm}
\end{equation}
The non-zero terms in \eqref{pm} occur when $j=i-1,i,i+1$.  The
$j=i$ term yields 
\begin{align*}
&\sum_{i=1}^n|\sigma_i|\int_{a_{i-1}}^{a_i}\int_{a_{i-1}-t}^{a_{i}
-t}p_m(x)\,dx\,dt\geq 
\sum_{i=1}^n|\sigma_i|\int_{a_{i-1}+\delta_m}^{a_i-\delta_m}\int_{-\delta_m}
^{\delta_m}p_m(x)\,dx\,dt\\
&=\sum_{i=1}^n|\sigma_i|\left(\frac{2\pi}{n}-2\delta_m\right)
\to \|\sigma\|_1 \text{ as } m\to\infty.
\end{align*}
When $j=i+1$, \eqref{pm} gives
\begin{align*}
&\left|\sum_{i=1}^{n-1}\epsilon_i\sigma_{i+1}\int_{a_{i}}^{a_{i+1}}
\int_{a_{i-1}-t}^{a_{i}-t}p_m(x)\,dx\,dt\right|\\
&=
\left|\sum_{i=1}^{n-1}\epsilon_i\sigma_{i+1}\int_{a_{i}}^{a_{i}+\delta_m}
\int_{-\delta_m}^{a_{i}-t}p_m(x)\,dx\,dt\right|
\leq n\delta_m\|\sigma\|_1/(2\pi)\to 0 \text{ as } m\to\infty.
\end{align*}
Similarly with the $j=i-1$ term in \eqref{pm}.
To prove (f), consider a sequence
$\{g_n\}\subset\bvt$ such that $\|g_k-g\|_1\to 0$ as $k\to\infty$.  
From  (g) in
Theorem~\ref{bvconv} we have $\widehat{f\ast g_k}(n)=\fhatn\hat{g_k}(n)$.
But $\{\widehat{g_k}\}$ converges to $\hat{g}$ as $k\to\infty$, 
uniformly on $\Z$
(\cite[I~Corollary~1.5]{katznelson}) so we can take the limit
$k\to\infty$ to complete the proof.
\qed

To see that (b) of Definition~\ref{abdefn} makes sense, take
a sequence $\{f_k\}\subset L^1(\T)$
such that $\|f_k-f\|_\T\to 0$.  The estimate
$\|f_k\ast g - f_l\ast g\|_\T\leq \|f_k-f_l\|_\T\|g\|_1$ from
Theorem~\ref{bvconv}(i) shows
$\{f_k\ast g\}$ converges to a unique element of $\alext$.
It is easy to see that this does not depend on the choice of
sequence $\{f_k\}$.
To see that (a) and (b)  agree, take $\{f_k\}$ and $\{g_k\}$ as in
Definition~\ref{abdefn}.
Then
\begin{eqnarray*}
\|f_k\ast g- f\ast g_k\|_\T & \leq & \|f_k\ast g-f_k\ast g_k\|_\T+
\|f\ast g_k-f_k\ast g_k\|_\T\\
 & \leq & \|f_k\|_\T\|g_k-g\|_1 + \|f_k-f\|_\T\|g_k\|_1.
\end{eqnarray*}
Since $\{\|f_k\|_\T\}$ and $\{\|g_k\|_1\}$ are bounded, letting
$k\to\infty$ shows that $f\ast g$ as defined by (a) and (b) are
the same.

The factorization theorem of Salem--Zygmund--Rudin--Cohen states that
if $E$ is any of the spaces $L^p(\T)$ for $1\leq p<\infty$ or
any of the spaces $C^k(\T)$ for $0\leq k<\infty$ then
$E=L^1(\T)\ast E$, i.e., for each $f\in E$ there exist $g\in L^1(\T)$ and
$h\in E$ such that $f=g\ast h$.  See \cite[7.5.1]{edwardsI}.
We have a similar result in $\alext$.
\begin{theorem}\label{factorization}
$\alext=L^1(\T)\ast \alext$.
\end{theorem}
\bigskip
\noindent
{\bf Proof:} Let $f\in\alext$.  Its primitive in $\balext$ is given by
$F(x)=\int_{-\pi}^xf$.  
Write $\ftilde =f - F(\pi)/(2\pi)$ and $\Ftilde(x)=\int_{-\pi}^x\ftilde$.
Then $\ftilde\in \alext$ and $\Ftilde\in C^0(\T)$ since $\intpi\ftilde=0$.
As $C^0(\T)=L^1(\T)\ast C^0(\T)$ there exist $g\in L^1(\T)$ and
$H\in C^0(\T)$ such that $\Ftilde=g\ast H$.
Differentiating both sides (\cite[Theorem~12]{talvilaconv}) 
gives $\ftilde=g\ast H'$.  Now let $c_1$ and $c_2$ be constants.
Then $(g+c_1)\ast (H'+c_2)=f-F(\pi)/(2\pi)+c_2\intpi g +2\pi c_1c_2$.
Let $c_2=1/(2\pi)$ and $c_1=(F(\pi)-\intpi g)/(2\pi)$ to complete the proof.  
\qed

This
theorem also follows from Theorem~22 and Note~25a in
\cite[\S32]{hewittrossII} since the approximate unit for $L^1(\T)$ is
also an approximate unit for $\alext$.
This is a sequence $\{k_n\}\subset L^1(\T)$ such that $\|k_n\|_1<M$ and
$\|f\ast k_n-f\|_1\to 0$ as $n\to\infty$, for each $f\in L^1(\T)$.
See \cite[28.51]{hewittrossII} and
Theorem~\ref{normconv} below.  This connection was pointed out by an
anonymous referee.

\begin{example}
{\rm
Using the method of Definition~\ref{abdefn}, it  does not seem 
possible to define the convolution on $\alext\times\alext$.
The following example shows there is no $k\in\R$ such that
$\|f\ast g\|_\T\leq k\|f\|_\T\|g\|_\T$ for all $f,g\in\alext$.
Let $f(t)=t^{-3}\sin(t^{-4})$ for $t\in(0,\pi)$,
let $f(t)=0$ for $t\in[-\pi,0]$ and extend $f$ periodically.
The primitive is 
$$
F(x)=\int_{-\pi}^xf=\left\{\begin{array}{cl}
0, & -\pi\leq x\leq 0\\
\frac{x^2}{4}\cos(x^{-4})-\frac{1}{2}\int_0^x t\cos(t^{-4})\,dt,
& 0<x<\pi,
\end{array}
\right.
$$
extended outside $[-\pi,\pi)$ so that $F\in\balext$.  Let
$f_n(t)=t^{-3}\sin(t^{-4})\chi_{((n\pi)^{-1/4},\pi)}(t)$. Extend
periodically outside $[-\pi,\pi)$ then $f_n\in \bvt$.
Now,
\begin{eqnarray*}
\|f-f_n\|_\T & = & \sup_{|\alpha-\beta|<2\pi}\left|\int_\alpha^\beta f-f_n\right| =
\max_{x,y\in[0,(n\pi)^{-1/4}]}|F(x)-F(y)|\\
& \leq & (n\pi)^{-1/2} \to 0 \text{ as } n\to\infty.
\end{eqnarray*}
Define $G(t)=t\sin(t^{-4})$ for $t\in[-\pi,0)$, $G(t)=0$ for
$t\in[0,\pi)$ and extend $G$ so that $G\in\balext$.  Let $g=G'$.
Then
$g\ast f_n(x)  =  \int_{(n\pi)^{-1/4}}^\pi g(x-t)f(t)\,dt$.
Using Theorem~\ref{fubini},
$$
\left|\int_0^\pi g\ast f_n(x)\,dx\right|  =  
\int_{(n\pi)^{-1/4}}^\pi \sin^2(t^{-4})\,\frac{dt}{t^2}
=\frac{1}{4}\int_{\pi^{-4}}^{n\pi}x^{-3/4}\sin^2(x)\,dx.
$$
Hence, $\|g\ast f_n\|_\T\to\infty$ as $n\to\infty$.
}
\end{example}

Since the convolution is linear in both arguments, associative over
$L^1(\T)$ and satisfies
the inequality $\|f\ast g\|_\T\leq \|f\|_\T \|g\|_1$, the convolution
maps
$\ast\fn\alext
\times L^1(\T)\to\alext$ and
$\alext$ is an $L^1(\T)$-module over the $L^1(\T)$ convolution algebra.
Trigonometric polynomials are dense in $\alext$ (Lemma~\ref{denselemma})
so $\alext$ is an essential Banach module.  
See \cite{daleslau} for the definitions.

A Segal algebra is a subalgebra of $L^1(\T)$ that is dense, translation
invariant, and continuous in norm.  See \cite{reiterstegeman}.  Since
$L^1(\T)$ is a subalgebra of $\alext$ the roles of the spaces are
reversed.  For
each $x\in\R$ and $f\in L^1(\T)$ we have $\|\tau_xf\|_1=\|f\|_1$ and
there is continuity
in $\|\cdot\|_1$.  Some properties of Segal algebras hold in this case.  For
example, if $f\in L^1(\T)$ then $\|f\|_\T\leq\|f\|_1$ 
(\cite[Proposition~6.2.3]{reiterstegeman}). 
However, it follows from Theorem~\ref{factorization} that
$L^1(\T)$ is not an ideal of $\alext$ 
(\cite[Proposition~6.2.4]{reiterstegeman}).

\section{Convergence}
The series $\sum_{-\infty}^\infty \fhat(n)e^{int}$ is known
as the Fourier series of $f$.  If $f$ is a  smooth enough
function then the Fourier series of $f$ converges to $f$.
There is  a substantial literature on pointwise convergence
of Fourier series.  For example, if the pointwise derivative
$f'(x)$ exists then the Fourier series converges to $f$  at $x$
\cite[Corollary~3.3.9]{grafakosI}.
It is a celebrated result of A.N.~Kolmogorov that there exists a
function $f\in L^1(\T)$ such that for each $t\in\T$ the sequence
$\sum_{-N}^N \fhat(n) e^{int}$ diverges as  $N\to\infty$ 
\cite[p.~80]{katznelson}.
L.~Carleson and R.A.~Hunt have proved that if $f\in L^p(\T)$ for some 
$1<p<\infty$ these symmetric
partial sums (given by convolution of $f$ with the Dirichlet kernel)
converge to $f$ almost everywhere.
For a proof see \cite[\S3.6]{grafakosI} together with
\cite[Chapter~11]{grafakosII}.  On the one-dimensional torus,  convergence 
of these symmetric partial 
sums 
to $f$ in the $p$-norm is equivalent to $L^p(\T)$ boundedness
of the conjugate function.  M.~Riesz has shown that the conjugate
function is bounded for $1<p<\infty$.  See \cite[\S3.5]{grafakosI}.
We will see below that these symmetric partial sums
converge to $f\in L^1(\T)$ in the Alexiewicz norm.  For
$f\in\alext$ we will show that the Fourier series
converges in the Alexiewicz norm with an appropriate summability factor.

First we consider summability kernels.

\begin{defn}\label{defnsummability}
A summability kernel is a sequence
$\{k_n\}\subset \bvt$ such that $\intpi k_n=1$,
$\lim_{n\to\infty}\int_{|s|>\delta}|k_n(s)|\,ds =0$ for each 
$0<\delta\leq\pi$ and
there is $M\in\R$ so that $\|k_n\|_1\leq M$ for all $n\in\N$.
\end{defn}
\begin{theorem}\label{normconv}
Let $f\in\alext$. 
Let $k_n$ be a summability kernel.
Then $\|f\ast k_n-f\|_\T\to 0$ as $n\to \infty$.
\end{theorem}
\bigskip
\noindent
{\bf Proof:}
Let $-\pi\leq \alpha<\beta\leq\pi$.  Then
\begin{align}
&\left|\int_{\alpha}^\beta \left[f\ast k_n(t)-f(t)\right]dt\right|
  =  \left|\int_{\alpha}^\beta \left[\intpi k_n(s)f(t-s)\,ds -f(t)
\intpi k_n(s) ds\right]dt\right|\notag\\
&=\left|\intpi k_n(s)\int_{\alpha}^\beta\left[f(t-s)-f(t)\right]
dt\,ds\right|\label{fubiniapproxid}\\
&\leq  \sup_{|s|<\delta}\|f-\tau_sf\|_\T  \int_{|s|<\delta}
|k_n(s)|\,ds +2\|f\|_\T\int_{\delta<|s|<\pi}
|k_n(s)|\,ds.\notag
\end{align}
The interchange of integrals in \eqref{fubiniapproxid} is accomplished using Theorem~\ref{fubini}
in the Appendix.
Due to continuity in the Alexiewicz norm,
given $\epsilon$, we can take $0<\delta<\pi$ small enough so that
$\|f-\tau_sf\|_\T<\epsilon$ for all $|s|<\delta$.  Hence,
$\|f\ast k_n-f\|_\T< M\epsilon + 2\|f\|_\T\int_{\delta<|s|<\pi}
|k_n(s)|\,ds$.  Letting $n\to\infty$ completes the proof. \qed

A commonly used summability kernel is the Fej\'er kernel,
$$
k_n(t)=\frac{1}{2\pi}
\sum_{k=-n}^n\left(1-\frac{|k|}{n+1}\right) e^{ikt}=\frac{1}
{2\pi(n+1)}\left[\frac{\sin((n+1)t/2)}{\sin(t/2)}\right]^2.
$$
See \cite{katznelson} for this and other summability kernels.
The classical summability kernels (de la Vall\'ee Poussin,
Poisson, Jackson) all satisfy the conditions
of Theorem~\ref{normconv}, which differ from Lebesgue integral
conditions by requiring the kernels be of bounded variation.
A sequence (or net) of functions satisfying the conclusion of
Theorem~\ref{normconv} is also called an approximate unit when
its Fourier series consists of a finite number of terms.
See \cite{reiterstegeman}.  The approximate units for $L^1(\T)$
are then approximate units for $\alext$.

\begin{lemma}\label{denselemma}
Let $f\in\alext$.  Then $f\ast e_n(x)=\fhatn e^{inx}$.
Let $g(t)=\sum_{-n}^n a_ke_k(t)$ for a sequence $\{a_k\}\subset\R$.
Then $f\ast g(x)=\sum_{-n}^n a_k\hat{f}(k) e^{ikx}$.
\end{lemma}
The proof follows from the identity $e_n(x-t)=e_n(x)e_n(-t)$ and
linearity of the integral.

The lemma allows us to prove that trigonometric polynomials are
dense in $\alext$ and gives a uniqueness result.
Let $k_n$ be the Fej\'er kernel and define $\sigma_n[f]=k_n\ast f$.
From Theorem~\ref{normconv} we have $\sigma_n[f]\to f$ in 
the Alexiewicz norm.  The Lemma shows $\sigma_n[f]$ is a trigonometric
polynomial.  Hence, the trigonometric polynomials are dense in $\alext$.
\begin{theorem}\label{uniqueness}
Let $f\in\alext$.
The
trigonometric polynomials are dense in $\alext$;
\begin{equation}
\sigma_n[f](t)=\frac{1}{2\pi}\sum_{k=-n}^n
\left(1-\frac{|k|}{n+1}\right)\hat{f}(k) e^{ikt}
\text{ and } \lim_{n\to\infty}\|f-\sigma_n[f]\|_\T=0.\label{fejerconv}
\end{equation}
If $\fhatn=0$ for all $n\in\Z$ then $f=0$.
\end{theorem}

Define the space of doubly indexed sequences converging to $0$ by 
$c_0=\{\sigma\fn\Z\to\R\mid \sigma_n=o(1)\text{ as } |n|\to \infty\}$.
Then $c_0$ is a Banach space under the uniform norm.
Distributions whose sequence of Fourier coefficients are in $c_0$ are known as
pseudo-functions.  Let 
$d=\{\sigma\fn\Z\to\R\mid \sigma_n=o(n)\text{ as } |n|\to \infty\}$.
Then $d$ is a Banach space under the norm $\|\sigma\|_d=\sup_{n\in\Z}
|\sigma_n|/(|n|+1)$.  In fact, $c_0$ and $d$ are  isometrically isomorphic, 
a linear
isometry being  given by $\sigma_n\mapsto\sigma_n/(|n|+1)$.
Note that a corollary to  Theorem~\ref{coefficient convergence} is
that if $\|f_j-f\|_\T\to 0$ then $\|\hat{f_j}-\fhat\|_d\to 0$.  The
following theorem summarizes the  properties of the mapping
$f\mapsto \fhat$ for  $f\in\alext$.
\begin{theorem}
Define ${\cal F}\fn\alext\to d$ by ${\cal F}[f]=\fhat$.  Then ${\cal F}$
is a bounded  linear transformation that  is  injective  but not 
surjective.
\end{theorem}
\bigskip
\noindent
{\bf Proof:} Linearity is given in Theorem~\ref{coefficients}(a), (b).
Part (f) of  the same theorem shows ${\cal  F}$ is bounded.
The uniqueness theorem (Theorem~\ref{uniqueness}) 
shows ${\cal F}$ is an injection.  If ${\cal F}$ were also a surjection
then  a consequence of the Open Mapping  Theorem is that 
there is $\delta>0$ such that $\|\fhat\|_d\geq\delta \|f\|_\T$
for all $f\in\alext$.
See \cite[Theorem~5.10]{rudin}.
For each $\alpha\in\R$ let $f_\alpha(t)=|t|^{-\alpha}{\rm sgn}(t)$ on 
$[-\pi,\pi)$.
Then $f_\alpha\in\alext$ if  and  only if  $\alpha<1$.
We  have $\hat{f_\alpha}(0)=0$ and for $n\not=0$ we get
$|\hat{f_\alpha}(n)|\leq 
2|n|\int_0^\pi t^{1-\alpha}\,dt =
2|n|\pi^{2-\alpha}/(2-\alpha)$.
Then  $\|\hat{f_\alpha}\|_d\leq 2\pi^{2-\alpha}/(2-\alpha) \to 2\pi$
as $\alpha\to 1^-$.
And, $\|f_\alpha\|_\T=\pi^{1-\alpha}/(1-\alpha)\to\infty$ as $\alpha\to
1^-$.  Hence, ${\cal F}$ cannot be  surjective. \qed

For $L^1$ Fourier series, $\fhat(n)=o(1)$ but the transformation
$f\mapsto \fhat$
is not onto $c_0$.  See \cite[Theorem~5.15]{rudin}.  

For $n\geq 0$ define the Dirichlet kernel
$D_n(t)=\sum_{-n}^n e^{ikt}=\sin[(n+1/2)t]/\sin(t/2)$.  Notice
that according to the definition in Theorem~\ref{normconv},
$D_n$ is not a summability kernel.  In fact, $\|D_n\|_1\sim(4/\pi^2)\log(n)$
as $n\to\infty$.  See \cite[p.~71]{katznelson}.  However, $\|D_n\|_\T$
are bounded.  This shows that  $D_n\ast f$ converges to $f$ in $\|\cdot\|_\T$
for $f\in L^1(\T)$.
\begin{theorem}\label{dirichlet}
The sequence $\|D_n\|_\T$ is bounded.
Let $f\in L^1(\T)$.  Then 
$\|D_n\ast f-f\|_\T\to 0$ as $n\to\infty$.
\end{theorem}
\bigskip
\noindent
{\bf Proof:}
Fix $n\in\N$.  
Since the function $t\mapsto \sin(t/2)$ is increasing on $[0,\pi]$
we have $\left|\int_{k\pi/(n+1/2)}^{(k+1)\pi/(n+1/2)}D_n(t)\,dt\right|\geq
\left|\int_{(k+1)\pi/(n+1/2)}^{(k+2)\pi/(n+1/2)}D_n(t)\,dt\right|$
for each integer $k\geq 0$.
We then have 
\begin{eqnarray*}
\|D_n\|_\T & = & 2\int_0^{2\pi/(2n+1)}\frac{\sin[(n+1/2)t]}{\sin(t/2)}dt
=2\int_0^{2\pi/(2n+1)}\sum_{k=-n}^n e^{ikt}\,dt\\
 & = & \frac{4\pi}{2n+1} +4\sum_{k=1}^n\frac{1}{k}\sin\left(\frac{
2\pi k}{2n+1}\right)\\
 & \leq & \frac{4\pi}{2n+1} + \frac{8\pi n}{2n+1} =4\pi.
\end{eqnarray*}
Let $f\in L^1(\T)$  and  let  $\epsilon>0$.
There is a trigonometric polynomial $p$ such that $\|f-p\|_1<\epsilon/(4\pi+1)$.
Let $n$  be greater than the degree of $p$.  Then from
Lemma~\ref{denselemma}  we  have
$D_n\ast p=p$. Using the estimate in Theorem~\ref{L1conv}(b)
\begin{eqnarray*}
\|D_n\ast f-f\|_\T & = & \|D_n\ast (f-p) +p-f\|_\T\\
 & \leq & \|D_n\|_\T\|f-p\|_1 + \|f-p\|_\T\\
 & \leq & (4\pi+1)\|f-p\|_1 <\epsilon.\qed
\end{eqnarray*}
\begin{example}\label{Dnexample}
{\rm
Since the Dirichlet kernels are not uniformly bounded in the $L^1$ norm
there is a function $f\in\alext$ such that $\|D_n\ast f-f\|_\T\not\to 0$.
To see this, for each $n\in\N$ define
$$
F_n(t) =\left\{\begin{array}{cl}
0, & 0\leq t\leq \pi\\
-\sin[(n+1/2)t], & -n\pi/(n+1/2)\leq t\leq 0\\
0, & -\pi\leq t\leq -n\pi/(n+1/2),
\end{array}
\right.
$$
with $F_n$ extended periodically.  Then $F_n\in\balext$ so $F_n'\in\alext$.

We have $D_n\ast F'_n=(D_n\ast F_n)'$ 
\cite[Proposition~4.2]{talvilaconv}.  Therefore, $\|D_n\ast F_n'\|_\T
=\max_{x,y\in[-\pi,\pi]}|D_n\ast F_n(y)-D_n\ast F_n(x)|$.  Note that
\begin{eqnarray*}
D_n\ast F_n(0) & = & \frac{2}{2n+1}\int_0^{n\pi}\frac{\sin^2(t)}
{\sin[t/(2n+1)]}\,dt
\geq 2\int_\pi^{n\pi}\sin^2(t)\,\frac{dt}{t}\\
 & = & \log(n)-\int_\pi^{n\pi}\cos(2t)\,\frac{dt}{t}.
\end{eqnarray*}
Hence, $D_n\ast F_n(0)\geq 0.5\log(n)$ for large enough $n$.
As well, 
$$
D_n\ast F_n(\pi)=\frac{(-1)^n}{2n+1}\int_0^{n\pi}\frac{\sin(2t)}{
\cos[t/(2n+1)]}\,dt.
$$
Since the function $t\mapsto \sec[t/(2n+1)]$ is positive and
increasing on $[0,n\pi]$
we have $\int_0^{n\pi}\sin(2t)
\sec[t/(2n+1)]\,dt<0$.

We now have
$\max_{x,y\in[-\pi,\pi]}|D_{2n}\ast F_{2n}(y)-D_{2n}\ast F_{2n}(x)|
\geq 0.5\log(n)$ for large enough $n$.  And, 
$\|F'_{2n}\|_\T=2$.  Hence, $\|D_n\|=\sup_{\|f\|_\T=1}\|D_n\ast f\|_\T$
are not uniformly bounded.  By the Uniform Boundedness Principle, there
exists $f\in\alext$ such that $\|D_n\ast f\|_\T$ is not bounded as $n\to
\infty$.  Therefore, $\|D_n\ast f-f\|_\T\not\to 0$.

Note  that if $f\in L^1(\T)$ then $\|D_n\ast f-f\|_1$ need  not  tend
to zero.  See \cite[p.~68]{katznelson}.
}
\end{example}

If $f\in L^1(\T)$ and $g\in L^\infty(\T)$ then Fej\'er's lemma states that
$\intpi f(t)g(nt)\,dt$ has limit $\fhat(0)\hat{g}(0)/(2\pi)$ as $n\to\infty$.
Since the multipliers for $\alext$ are  the functions of bounded variation
we have the
following version  for the continuous primitive integral.
\begin{prop}\label{fejerlemma}
Let $f\in\alext$ and $g\in\bvt$.  Then $\intpi f(t)g(nt)\,dt= 
o(n)$
as $n\to\infty$.  The order estimate is sharp.
\end{prop}
\bigskip
\noindent
{\bf Proof:}
The trigonometric polynomials are dense in $\alext$ (Lemma~
\ref{denselemma}) so there are sequences of trigonometric polynomials
$\{p_\ell\}$ and $\{q_m\}$ such that $\|f-p_\ell\|_\T\to 0$
and $\|g-q_m\|_\T\to0$ as $\ell, m\to\infty$.  Write
\begin{eqnarray*}
\intpi f(t)g(nt)\,dt & = & \intpi[f(t)-p_\ell(t)]g(nt)\,dt+\intpi p_\ell(t)
[g(nt)-q_m(nt)]\,dt\\
 & & \quad +\intpi p_\ell(t)q_m(nt)\,dt.
\end{eqnarray*}
Use the H\"older inequality, Proposition~\ref{holder}.
Then $|\intpi[f(t)-p_\ell(t)]g(nt)\,dt|\leq n\|f-p_\ell\|_\T \|g\|_\bv$.
Take $\ell\in\N$ large enough so that $\|f-p_\ell\|_\T \|g\|_\bv<\epsilon$.
And, $|\intpi p_\ell(t)
[g(nt)-q_m(nt)]\,dt|\leq \|p_\ell\|_\bv \|g-q_m\|_\T$.
Take $m\in\N$ large enough so that $\|p_\ell\|_\bv \|g-q_m\|_\T<\epsilon$.
Now write $q_m=\sum_{-m}^m a_k e_{-k}$.  Then
$\intpi p_\ell(t)q_m(nt)\,dt=\sum_{-m}^m a_k\hat{p}(nk)$.
By the Riemann--Lebesgue lemma this tends to $a_0\hat{p}(0)$ as $n\to\infty$.
Hence, $\intpi f(t)g(nt)\,dt= o(n)$.
The estimate is sharp by Theorem~\ref{coefficients}(g).
\qed

Since the topological dual of $\alext$ is $\bvt$ we can view
functions of bounded variation as continuous linear functionals
on $\alext$.  For $f\in\alext$ and $g\in\bvt$ we define
a linear  functional $g\fn\alext\to\R$ by $g[f]=\intpi fg$.  
If $\{f_n\}\subset\alext$ such that $\|f_n-f\|_\T\to 0$
then by the H\"older inequality 
$$
|g[f_n]-g[f]|=\left|\intpi (f_n-f)g\right|\leq \|f_n-f\|_\T\|g\|_\bv\to0.
$$
Hence $g$ is a continuous linear functional on $\alext$. 

The following Parseval equality states that for every $g\in\bv$ we have
$\sigma_n[g]\to g$ in  the weak$^\ast$ topology.
\begin{theorem}\label{parseval}
If $f\in\alext$ and $g\in\bvt$ then
$$
g[f]=\langle f,g\rangle=\intpi fg=
\lim_{n\to\infty}\sum_{k=-n}^n\left(1-\frac{|k|}{n+1}\right)
\hat{f}(k)\hat{g}(k).
$$
\end{theorem}
The proof is essentially the same as the version for Lebesgue integrals
given in \cite[p.~37]{katznelson}.  Note that for $f\in\alext$ and
$g\in\bvt$, the series $\sum_{-\infty}^\infty
\hat{f}(k)\hat{g}(k)$ need not converge.  This follows from the
sharp growth estimates $\fhat(k)=o(k)$ (Theorem~\ref{coefficients}(g))
and $\hat{g}(k)=O(1/k)$ (\cite[I~Theorem~4.5]{katznelson}).  
To see this, let $g(t)=\sqrt{1-(t/\pi)^2}$
on $(0,\pi)$ and extend $g$ as an odd periodic function.  Then
$$
\hat{g}(k)=\frac{2i}{k} - \frac{2i}{\pi^2 k}\int_0^\pi\frac{t\cos(kt)\,dt}
{\sqrt{1-(\pi/t)^2}}.
$$
By the Riemann--Lebesgue lemma, $\hat{g}(k)\sim 2i/k$ as $k\to\infty$.

The Parseval equality lets us characterize sequences of Fourier coefficients
of functions  in $\bvt$.
\begin{theorem}\label{bvcharacterization}
Let $\{a_n\}_{n\in\Z}$ be a sequence in $\C$.  The following are
equivalent: (a)  There exists $g\in\bvt$ and  $c\geq 0$ such that
$\|g\|_\bv\leq c$ and $\hat{g}(n)=a_n$ for each $n\in\Z$.  (b) For
all trigonometric polynomials $p$ we have $|\sum_{-\infty}^\infty
\hat{p}(n) \overline{a_n}|\leq c\|p\|_\T$.
\end{theorem}
\begin{corollary}\label{bvseries}
A trigonometric series $S(t)\sim \sum a_n e^{int}$ is the Fourier series
of some function $g\in\bvt$ with $\|g\|_\bv\leq c$ if and only if
$\|\sigma_n[S]\|_\bv\leq c$ for all $n\in\Z$.
\end{corollary}
The proof is essentially the same as that for Theorem~7.3 in 
\cite[p.~39]{katznelson}.

\section{Appendix}\label{appendix}
The following type of Fubini theorem generalizes a similar one
for the Henstock--Kurzweil and wide Denjoy integral in 
\cite[Theorem~58]{celidze}.
\begin{theorem}\label{fubini}
Let $f\in\alext$.  Let $g\in\bvt$.  If $-\infty<a<b<\infty$ then
$
\int_a^b\intpi f(x-y)g(y)\,dy\,dx =
\intpi\int_a^b f(x-y)g(y)\,dx\,dy.
$
\end{theorem}
\bigskip
\noindent
{\bf Proof:}
Let $F\in\balext$ be the primitive of $f$.
Integrating by parts gives $\intpi f(x-y)g(y)\,dy =F(x+\pi)g(-\pi)
-F(x-\pi)g(\pi)+\intpi F(x-y)\,dg(y)$.  Now use the periodicity of
$g$ to write 
$$
\int_a^b\intpi f(x-y)g(y)\,dy\,dx =\left(\int_{a+\pi}^{b+\pi}F
-\int_{a-\pi}^{b-\pi}F\right)g(\pi)+\int_a^b\intpi F(x-y)\,dg(y)\,dx.
$$
A linear change of variables and integration by parts gives
\begin{align*}
&\intpi\int_a^b f(x-y)g(y)\,dx\,dy  =  \intpi g(y)\int_{a-y}^{b-y}f(x)
\,dx\,dy\\
& =\intpi \left[F(b-y)-F(a-y)\right]g(y)\,dy\\
& =\intpi[F(b-y)-F(a-y)]\,dy\,g(\pi) - \intpi\int_{-\pi}^y[F(b-z)
-F(a-z)]\,dz\,dg(y)\\
& =\left(\int_{b-\pi}^{b+\pi}F-\int_{a-\pi}^{a+\pi}F\right)g(\pi)
 + \intpi\left(\int_a^bF(x-y)\,dx-\int_{a+\pi}^{b+\pi}F\right)dg(y)\\
& =\left(
\int_{a+\pi}^{b+\pi}F-\int_{a-\pi}^{b-\pi}F\right)g(\pi) +
\intpi\int_{a}^{b}F(x-y)\,dx\,dg(y) - \int_{a+\pi}^{b+\pi}F\intpi dg.
\end{align*}
The usual Fubini theorem gives $\intpi\int_{a}^{b}F(x-y)\,dx\,dg(y)
=\int_{a}^{b}\intpi F(x-y)\,dg(y)\,dx$ since
$|\intpi\int_{a}^{b}F(x-y)\,dx\,dg(y)|\leq \max_{x\in[a-\pi, b+\pi]}|F(x)| Vg$.
As $g$ is periodic, $\intpi dg=0$.\qed

\end{document}